\documentclass[12pt]{amsart}
\usepackage{amsfonts,amssymb,color}
\usepackage[mathscr]{eucal}
\usepackage{amsmath, amsthm}
\usepackage{mathrsfs}
\usepackage{amsbsy}
\usepackage{dsfont}
\usepackage{bbm}
\usepackage{wasysym}
\usepackage{url}
\usepackage[active]{srcltx} 
\input xypic
\xyoption{all}
\begin{document}
\baselineskip = 5mm
\newcommand \ZZ {{\mathbb Z}} 
\newcommand \NN {{\mathbb N}} 
\newcommand \QQ {{\mathbb Q}} 
\newcommand \RR {{\mathbb R}} 
\newcommand \CC {{\mathbb C}} 
\newcommand \PR {{\mathbb P}} 
\newcommand \AF {{\mathbb A}} 
\newcommand \uno {{\mathbbm 1}}
\newcommand \Le {{\mathbbm L}}
\newcommand \bcA {{\mathscr A}}
\newcommand \bcB {{\mathscr B}}
\newcommand \bcC {{\mathscr C}}
\newcommand \bcD {{\mathscr D}}
\newcommand \bcE {{\mathscr E}}
\newcommand \bcF {{\mathscr F}}
\newcommand \bcI {{\mathscr I}}
\newcommand \bcM {{\mathscr M}}
\newcommand \bcS {{\mathscr S}}
\newcommand \bcT {{\mathscr T}}
\newcommand \bcU {{\mathscr U}}
\newcommand \bcX {{\mathscr X}}
\newcommand \bcY {{\mathscr Y}}
\newcommand \Spec {{\rm {Spec}}}
\newcommand \Pic {{\rm {Pic}}}
\newcommand \Alb {{\rm {Alb}}}
\newcommand \Corr {{Corr}}
\newcommand \Sym {{\rm {Sym}}}
\newcommand \LSym {L{\rm {Sym}}}
\newcommand \Alt {{\rm {Alt}}}
\newcommand \cha {{\rm {char}}}
\newcommand \tr {{\rm {tr}}}
\newcommand \im {{\rm im}}
\newcommand \Hom {{\rm Hom}}
\newcommand \Aut {{\rm Aut}}
\newcommand \bHom {{\bf Hom}}
\newcommand \iHom {{\underline {\rm Hom}}}
\newcommand \colim {{{\rm colim}\, }} 
\newcommand \hocolim {{{\rm hocolim}\, }} 
\newcommand \Hocolim {{{\rm Hocolim}\, }} 
\newcommand \coeq {{{\rm coeq}\, }} 
\newcommand \eq {{{\rm eq}\, }} 
\newcommand \End {{\rm {End}}}
\newcommand \coker {{\rm {coker}}}
\newcommand \id {{\rm {id}}}
\newcommand \Id {{\rm {Id}}}
\newcommand \Ob {{\rm Ob}}
\newcommand \Cyl {{\rm Cyl\; }}
\newcommand \Map {{\rm {Map}}}
\newcommand \map {{\rm {map}}}
\newcommand \op {{^{\rm op}}} 
\newcommand \cor {{\rm {cor}}}
\newcommand \res {{\rm {res}}}
\newcommand \cone {{\rm {cone}}}
\newcommand \Spt {{\rm Spt}}
\newcommand \Ev {{\rm Ev}}
\newcommand \Com {{\rm Com}}
\newcommand \Tot {{\rm Tot}}
\newcommand \dom {{\rm dom}}
\newcommand \cell {{\rm cell}}
\newcommand \inj {{\rm inj}}
\newcommand \cof {{\rm cof}}
\newcommand \proj {{\rm proj}}
\newcommand \codom {{\rm codom}}
\newcommand \Ord {{\it Ord}}
\newcommand \mg {{\mathfrak m}}
\newcommand \sg {{\Sigma }}
\newcommand \CHM {{\mathscr C\! \mathscr M}}
\newcommand \DM {{\mathscr D\! \mathscr M}}
\newcommand \MM {{\mathscr M\! \mathscr M}}
\newcommand \pt {{*}} 
\newcommand \cf {{\rm cf}} 
\newcommand \SH {{SH}} 
\newcommand \MSH {{\bf SH}} 
\newcommand \de {{\vartriangle }} 
\newcommand \deop {{\vartriangle }^{op}} 
\newcommand \SSets {{\deop \mathscr S\! ets}} 
\newcommand \Sets {{\mathscr S\! ets}} 
\newcommand \Sm {{\mathscr S\! m}} 
\newcommand \Sch {{\mathscr S\! ch}} 
\newcommand \Gm {{{\mathbb G}_{\rm m}}}
\newcommand \cX {{\mathcal X}} 
\newcommand \cY {{\mathcal Y}} 
\def\blue {\color{blue}}
\def\red {\color{red}}
\newtheorem{theorem}{Theorem}
\newtheorem{lemma}[theorem]{Lemma}
\newtheorem{corollary}[theorem]{Corollary}
\newtheorem{proposition}[theorem]{Proposition}
\newtheorem{remark}[theorem]{Remark}
\newtheorem{definition}[theorem]{Definition}
\newtheorem{conjecture}[theorem]{Conjecture}
\newtheorem{example}[theorem]{Example}
\newcommand \lra {\longrightarrow}
\newenvironment{pf}{\par\noindent{\em Proof}.}{\hfill\framebox(6,6) \par\medskip}
\title[Positive model structures]
{\bf Positive model structures for abstract symmetric spectra}
\author{S. Gorchinskiy, V. Guletski\u \i }

\date{11 May 2015}

\begin{abstract}
\noindent We give a general method of constructing positive stable model structures for symmetric spectra over an abstract simplicial symmetric monoidal model category. The method is based on systematic localization, in Hirschhorn's sense, of a ceratin positive projective model structure on spectra, where positivity basically means the truncation of the zero slice. The localization above is by the set of stabilizing morphisms, or their truncated version.
\end{abstract}









\maketitle

\section{Introduction}
\label{s-intro}

The aim of this paper is to give a systematic account of the method of constructing positive model structures for abstract symmetric spectra, used to prove one of the key theorems in \cite{GG2}. Let first $\bcS $ be the category of topological symmetric spectra in the sense of \cite{HSS}, and let $\bcT $ be the homotopy category of $\bcS $ with respect to the stable model structure in it. Then $\bcT $ is equivalent to the standard topological stable homotopy category, whose Hom-groups encode the stable homotopy groups of $CW$-complexes. As it was shown in \cite{EKMM} (see also \cite{MMSS}), the category $\bcS $ admits another one, so-called positive, model structure whose homotopy category is the same as $\bcT $, but the positivity of this new structure gives rise to many good properties missing in the standard stable model structure. For example, if $X$ is a topological symmetric spectrum, which is cofibrant in the positive model structure, then the natural morphism from the $n$-th homotopy symmetric power of $X$ onto the honest $n$-th symmetric power of $X$ is a stable weak equivalence, loc.cit. The latter result is important for our understanding of the stable homotopy groups through the Barratt-Priddy-Quillen theorem, see the modern approach in \cite{Schlichtkrull}. Another essential application of positive model structures in topology is that it yields a convenient model structure for commutative ring spectra, see \cite{Shipley}.

On the other hand, following \cite{Hovey2}, one can get a general method for constructing stable homotopy categories, equally appropriate in topology and in $\AF ^1$-homotopy theory, where the initial category $\bcC $ is nothing but the category of simplicial Nisnevich sheaves on smooth schemes over a base, see \cite{Jardine}. We start with a closed symmetric monoidal model category $\bcC $, which is, in addition, left proper and cellular, then take a cofibrant object $T$ in $\bcC $ and look at the category of symmetric $T$-spectra $\bcS $ over $\bcC $. This category $\bcS $ possesses a stable model structure, and the corresponding homotopy category $\bcT $ generalizes the topological stable homotopy category and the motivic one, loc.cit. A natural question is then how to extend the method of constructing positive model structures developed in topology to the level of generality, high enough to be applicable in motivic algebraic geometry, and in other reasonable settings.

In the paper, we give an affirmative answer to this question and show a universal method of constructing many positive structures, adjustable to particular needs. Basically, we follow the method in \cite{EKMM} and \cite{MMSS}, keeping the level of generality as high as possible. A new thing, however, is that we systematically exploit the technique of localization of model categories from \cite{Hirsch}, which allows us to make the approach to be more conceptual and put an order on various model structures naturally arising in our considerations. In nutshell, we first take a projective model structure, truncate it in its $0$-slice, or any finite number of slices starting from the zero one, and then localize the truncated model structure by the stabilizing Hovey's $\zeta $-morphisms between appropriately shifted $T$-spectra.

The application of positive model structures in \cite{GG2} goes as follows. Let $X$ be an object cofibrant with respect to the positive projective model structure in $\bcS $. Then the natural morphism from the $n$-th homotopy symmetric power of $X$ to its honest $n$-th symmetric power is a stable weak equivalence of symmetric spectra. As a consequence, symmetric powers preserve stable weak equivalences between positively cofibrant objects in $\bcS $. This result generalizes Lemma 15.5 in \cite{MMSS}, and allows us to derive symmetric powers in the abstract stable homotopy category $\bcT $, see \cite{GG2}. The level of generality is high enough to apply the result in the Morel-Voevodsky stable category of motivic symmetric spectra over a field, loc.cit. Positive model structures were utilized in \cite{Hornbostel} in the context of $\AF ^1$-homotopy theory of schemes. They are also needed to compare the geometric symmetric powers of motivic spectra with their left derived symmetric powers, see \cite{Joe}. In \cite{CD} positive model structures were used for the study of commutative monoids in an abstract symmetric monoidal model category. In \cite{PavlovScholbach} the methods and results of the present work are extended to algebra spectra over symmetric operads.

The paper is organized as follows. In Section \ref{positive} we set up what exactly we want to construct, and fix notation and terminology. We recall some basic definitions on abstract symmetric spectra in Section \ref{positive}, but the reader is advised to use Hovey's article \cite{Hovey2} to repeat the details. In Section \ref{positiveproj} we present our concept of positive stable model structures as systematic localizations of positive projective model structures on symmetric spectra. We have chosen to start with projective model structure, but injective model structures are good for our purposes too. Section \ref{loops} is devoted to deducing the needed results on loop-spectra in the abstract setting. Finally, in Section \ref{stablepositive} we prove the main result (Theorem \ref{main}) saying that weak equivalences in the stable model structure are the same as weak equivalences in the positive model structure. This implies that the resulting stable homotopy category is the same.

\medskip

{\sc Acknowledgements.} The authors are grateful to Peter May, who has drawn our attention to positive model structures in topology, and to Joseph Ayoub for useful comments on homotopy types under the action of finite groups. The paper is written in the framework of the EPSRC grant EP/I034017/1. The first named author acknowledges the support of the grants MK-5215.2015.1, NSh-2998.2014.1, RFBR 13-01-12420, 14-01-00178, Dmitry Zimin's Dynasty Foundation and the subsidy granted to the HSE by the Government of the Russian Federation for the implementation of the Global Competitiveness Program.

\section{Positive model structures: what to construct?}
\label{positive}

First we need to explain what do we mean by an abstract stable homotopy category. Our viewpoint is that it should be understood as the homotopy category of the category of symmetric spectra over a given simplicial model monoidal category $\bcC $, stabilizing smashing with $T$, where $T$ is a cofibrant object $T$ in $\bcC $. Such a general gadget generalizes both the topological stable homotopy category and the motivic one due to Morel and Voevodsky. Nowadays, in both cases, we should work with symmetric spectra as they provide a set of powerful monoidal properties of spectra, useful in applications. In our considerations we depart from the paper \cite{Hovey2}, which is basic to us.

Let $\bcC $ be a closed symmetric monoidal model category with the monoidal product $\wedge $. This notation is the tradition coming from the pointed setting needed to make the homotopy category of spectra to be additive. Respectively, the coproduct will be denoted by $\vee $.

Next, we assume that the model structure in $\bcC $ is left proper and cellular. Left properness means that the push-out of a weak equivalence along a cofibration is again a weak equivalence, and cellularity means that $\bcC $ is cofibrantly generated by a set of generating cofibrations $I$ and a set of trivial generating cofibrations $J$, the domains and codomains of morphisms in $I$ are compact relative to $I$, the domains of morphisms in $J$ are small relative to the cofibrations, and cofibrations are effective monomorphisms. To avoid any misunderstanding in using this complicated terminology we would recommend the reader to consult with \cite{Hovey1}, \cite{Hovey2} and \cite{Hirsch}. Suppose, moreover, that the domains of the generating cofibrations $I$ in $\bcC $ are cofibrant, which is needed to satisfy the assumptions of Theorem 8.11 in \cite{Hovey2}.

For simplicity, we shall also assume that $\bcC $ is simplicial, and that the simplicial structure is compatible with the structure of a closed symmetric monoidal model category. This will be used in the proofs of Proposition \ref{lemma5} and Corollary \ref{kozerog} merely in order to avoid the bulky work with functional complexes. However, Proposition \ref{lemma5} and Corollary \ref{kozerog}, as well as the main Theorem \ref{main}, are true without this assumption.

Let $\sg $ be a disjoint union of symmetric groups $\sg _n$ for all $n\geq 0$, where $\sg _0=\emptyset $ and all groups are considered as one object categories. Let $\bcC ^{\sg }$ be the category of symmetric sequences over $\bcC $, i.e. functors from $\sg $ to $\bcC $. Since $\bcC $ is closed symmetric monoidal, so is the category $\bcC ^{\sg }$. The monoidal product in $\bcC ^{\sg }$ is given by the formula
  $$
  (X\wedge Y)_n=\vee _{i+j=n}\sg _n\times _{\sg _i\times \sg _j}(X_i\wedge Y_j)\; ,
  $$
where $\sg _n\times _{\sg _i\times \sg _j}(X_i\wedge Y_j)$ is nothing but $\cor ^{\sg _n}_{\sg _i\times \sg _j}(X_i\wedge Y_j)$ in terms of \cite{GG2}, and the action of $\sg _n$ is standard, see \cite{HSS} or \cite{Hovey2}.

Let $T$ be a cofibrant object in $\bcC $, and let $S(T)$ be the free monoid on the symmetric sequence $(\emptyset ,T,\emptyset ,\emptyset ,\dots )$, i.e. the symmetric sequence
  $$
  S(T)=(T^0,T^1,T^2,T^3,\dots )\; ,
  $$
where $T^0=\uno $ is the unit, $T^1=T$ and $\sg _n$ acts on $T^n$ by permutation of factors. The whole point is that the monoid $S(T)$ is commutative. Then a symmetric spectrum is nothing but a module over $S(T)$ in $\bcC ^{\sg }$. Explicitly, a symmetric spectrum $X$ is a sequence of objects
  $$
  X_0\; ,\; \; X_1\; ,\; \; X_2\; ,\; \; X_3\; ,\; \dots
  $$
in $\bcC $ together with $\sg _n$-equivariant morphisms
  $$
  X_n\wedge T\lra X_{n+1}\; ,
  $$
such that for all $n,i\geq 0$ the composite
  $$
  X_n\wedge T^i\lra X_{n+1}\wedge T^{i-1}\to \dots \to X_{n+i}
  $$
is $\sg _n\times \sg _i$-equivariant.

Let
  $$
  \bcS =\Spt ^{\sg }(\bcC ,T)
  $$
be the category of symmetric spectra over $\bcC $ stabilizing the functor
  $$
  -\wedge T:\bcC \lra \bcC \; .
  $$
There is a natural closed symmetric monoidal structure on $\bcS$ given by the product of modules over the commutative monoid $S(T)$.

A model structure on $\bcS $ can be constructed as a localization of the so-called projective model structure coming from the model structure on $\bcC $, using the main result of \cite{Hirsch}.

Namely, for any non-negative $n$ we consider the evaluation functor
  $$
  \Ev _n:\bcS \lra \bcC
  $$
sending any symmetric spectrum $X$ to its $n$-slice. Each $\Ev _n$ has a left adjoint
  $$
  F_n:\bcC \lra \bcS \; ,
  $$
which can be constructed as follows. Let $\tilde F_n$ be a functor sending any object $X$ in $\bcC $ to the symmetric sequence
  $$
  (\emptyset ,\dots ,\emptyset ,\sg _n\times X,\emptyset ,\emptyset ,\dots )\; ,
  $$
where $\emptyset $ is the initial object in $\bcC $. Then
  $$
  F_nX=\tilde F_nX\wedge S(T)\; ,
  $$
see \cite[Def.7.3]{Hovey2}.

Let now
  $$
  I_T=\cup _{n\geq 0}F_nI\hspace{8mm}\hbox{and}\hspace{8mm}
  J_T=\cup _{n\geq 0}F_nJ\; ,
  $$
where $F_nI$ is the set of all the morphisms of type $F_nf$, $f\in I$, and the same for $F_nJ$. Let also
  $$
  W_T
  $$
be the set of projective weak equivalences, i.e. level weak equivalences, which means that for any morphism $f:X\to Y$ in $W_T$ the morphism $f_n:X_n\to Y_n$ is a weak equivalence in $\bcC $ for all $n\geq 0$.

For technical reasons, we prefer to use different symbols to denote a category and a model structure in it. The projective model structure
  $$
  \bcM =(I_T,J_T,W_T)
  $$
is generated by the set of generating cofibrations $I_T$ and the set of trivial generating cofibrations $J_T$. As the model structure in $\bcC $ is left proper and cellular, the projective model structure in $\bcS $ is left proper and cellular too, \cite{Hovey2}. In particular, the class of cofibrations in $\bcM $ is equal to the class $I_T$-cof.

For any two non-negative integers $n$ and $m$, $m\geq n$, the group $\sg _{m-n}$ is canonically embedded into the group $\sg _n$, such that for any object $X$ in $\bcC $ it acts on $X\wedge T^{m-n}$ permuting factors in $T^{m-n}$. Then $F_nX$ can be computed by the formula
  $$
  (F_nX)_m=\cor ^{\sg _m}_{\sg _{m-n}}(X\wedge T^{m-n})\; ,
  $$
see \cite[\S 7]{Hovey2}. In particular,
  $$
  \Ev _{n+1}F_nX=\cor ^{\sg _{n+1}}_{\sg _1}(X\wedge T)=
  \sg _{n+1}\times (X\wedge T)\; .
  $$

Let now
  $$
  \zeta ^X_n:F_{n+1}(X\wedge T)\lra F_n(X)
  $$
be the adjoint to the morphism
  $$
  X\wedge T\lra \Ev _{n+1}F_nX=\sg _{n+1}\times (X\wedge T)
  $$
induced by the canonical embedding of $\sg _1$ into $\sg _{n+1}$.

For any set of morphisms $U$ let $\dom (U)$ and $\codom (U)$ be the set of domains and codomains of morphisms from $U$, respectively. Let then
  $$
  S=\{ \zeta ^X_n\; \mid \; X\in \dom (I)\cup \codom (I)\; ,\; \; n\geq 0\}
  $$
be the set of stabilizing morphisms. Then a {\it stable model structure}
  $$
  \bcM _S=(I_T,J_{T,S},W_{T,S})
  $$
in $\bcS $ is defined to be the Bousfield localization of the projective model structure with respect to the class $S$ in the sense of \cite{Hirsch}. It is generated by the same set of generating cofibrations $I_T$, and by a new set of trivial generating cofibrations $J_{T,S}$. Here $W_{T,S}$ is the set of stable weak equivalences, i.e. new weak equivalences obtained as a result of the localization.

Let
  $$
  \bcT =\bcS [W_{T,S}^{-1}]
  $$
be the localization of $\bcS $ with respect to the class $W_{T,S}$, i.e. the homotopy category of $\bcS $ with respect to weak equivalences in $W_{T,S}$. Then we call $\bcT $ to be an abstract stable homotopy category of symmetric spectra over $\bcC $ which stabilizes smashing by $T$. As the functor $(-\wedge T)$ is a Quillen autoequivalence of $\bcS $ with respect to the model structure $\bcM _S$, it induces an autoequivalence on the homotopy category $\bcT $, as required.

By Hovey's result, \cite{Hovey2}, the homotopy category $\bcT $ is equivalent to the homotopy category of ordinary $T$-spectra provided the cyclic permutation on $T\wedge T\wedge T$ is left homotopic to the identity morphism.

\bigskip

Let now
  $$
  S^+=\{ \zeta ^X_n\; \mid \; X\in \dom (I)\cup \codom (I)\; ,\; \; n>0\}
  $$
be the positive stabilizing set. Our aim is actually to find a new model structure $\bcM ^+$, generated by a new set $I_T^+$ of generating cofibrations, and a new set of generating trivial cofibrations $J_T^+$, having a new set of weak equivalences $W_T^+$
  $$
  \bcM ^+=(I_T^+,J_T^+,W_T^+)\; ,
  $$
such that weak equivalences in $\bcM ^+$ would be those morphisms $f:X\to Y$ in which $f_n:X_n\to Y_n$ is a weak equivalence in $\bcC $ for all $n>0$, and if
  $$
  \bcM ^+_{S^+}=(I^+_T,J^+_{T,S^+},W^+_{T,S^+})
  $$
is a localization of $\bcM ^+$ with respect to the above set $S^+$ then
  $$
  W_{T,S}=W^+_{T,S^+}\; .
  $$

Since now the desired model structure $\bcM ^+_{S^+}$ will be called a {\it positive stable model structure} whose fibrations, cofibrations and weak equivalences will be called {\it positive fibrations, cofibrations and stable weak equivalences}.

\section{Positive projective model structures}
\label{positiveproj}

Let
  $$
  I_T^+=\cup _{n>0}F_n(I)\; ,
  $$
  $$
  J_T^+=\cup _{n>0}F_n(J)
  $$
and let
  $$
  W_T^+
  $$
be the set of morphisms $f:X\to Y$, such that $f_n:X_n\to Y_n$ is a weak equivalence in $\bcC $ for all $n>0$. First we will prove a proposition saying that $I_T^+$, $J_T^+$ and $W_T^+$ generate a model structure in $\bcS $.

\begin{proposition}
\label{thA}
The above sets $I_T^+$, $J_T^+$ and $W^+_T$ do satisfy the conditions of Theorem 2.1.19 in \cite{Hovey1}, so that they generate a model structure, denoted by $\bcM ^+$ with the set of generating cofibrations $I_T^+$, the set of trivial generating cofibrations $J_T^+$, and whose weak equivalences are $W_T^+$. In particular, the set of cofibrations in $\bcM ^+$ is the set {\rm $I_T^+$-cof}\footnote{we will systematically use the terminology from $\S $2.1 of the book \cite{Hovey1}}, the set of trivial cofibrations is {\rm $J_T^+$-cof}, and weak equivalences in $\bcM ^+$ are $W_T^+$.
\end{proposition}

\begin{pf}
We will use the fact that $\bcM =(I_T,J_T,W_T)$, and so the sets $I_T$, $J_T$ and $W_T$ satisfy the conditions of Theorem 2.1.19 in \cite{Hovey1}.

\bigskip

\noindent {\it First condition}

\medskip

\noindent The first condition from Theorem 2.1.19 in \cite{Hovey1} is satisfied automatically.

\bigskip

\noindent {\it Second condition}

\medskip

\noindent Since
  $$
  I_T^+\subset I_T\; ,
  $$
we get
  $$
  \dom (I_T^+)\subset \dom (I_T)\; ,
  $$
and
  $$
  \hbox{$I_T^+$-cell} \subset \hbox{$I_T$-cell} \; .
  $$

By the property 2 from Hovey's theorem, applied to $\bcM $, we have that $\dom (I_T)$ are small relative to $I_T$-cell. Since $\dom (I_T^+)\subset \dom (I_T)$, even more so the set $\dom (I_T^+)$ is small relative to $I_T$-cell. As $I_T^+$-cell is a subset in $I_T$-cell, even more so the set $\dom (I_T^+)$ is small relative to the smaller class $I_T^+$-cell.

\bigskip

\noindent {\it Third condition}

\medskip

\noindent Everything is the same as in the case of the second condition, but we need to replace $I$ by $J$.

\bigskip

\noindent {\it Fourth condition}

\medskip

\noindent First we look at the chain of the obvious inclusions
  $$
  \hbox{$J_T^+$-cell}\subset \hbox{$J_T$-cell}\subset W_T\subset W_T^+\; .
  $$
Now we need to show that $\hbox{$J_T^+$-cell}\subset \hbox{$I_T^+$-cof}$. Notice that the class $J_T^+$-cell consists of transfinite compositions of push-outs of morphisms from $J_T^+$ and the class $I_T^+$-cof is closed under transfinite compositions and push-outs, see the proof of Lemma 2.1.10 on page 31 in \cite{Hovey1}. This is why, in order to show that $\hbox{$J_T^+$-cell}\subset \hbox{$I_T^+$-cof }$, it is enough to prove that $J_T^+\subset \hbox{$I_T^+$-cof }$.

\medskip

We need some more terminology. Let $\bcX $ be a category, and let $A$ and $B$ be two classes of morphisms in it. We will say that the the pair $\{ A,B\} $ has the lifting property (LP, for short) if for any morphism $f:X\to Y$ from $A$, and any morphism $g:U\to V$ from $B$, and any commutative square
  $$
  \diagram
  X \ar[dd]_-{f} \ar[rr]^-{} & & U \ar[dd]^-{g} \\ \\
  Y \ar@{.>}[uurr]^-{\exists \gamma }\ar[rr]^-{} & & V
  \enddiagram
  $$
there exists a morphism $\gamma $ keeping the diagram to be commutative.

Let now $\bcX $ and $\bcY $ be two categories, and let $F:\bcX \rightleftarrows \bcY :G$ be
two adjoint functors, $F$ from the left, and $G$ from the right. Let $A$ be a class of morphisms in $\bcX $, and let $B$ be a class of morphisms in $\bcY $. Then $\{ A,G(B)\} $ has LP if and only if $\{ F(A),B\} $ has LP.

Using this, and also taking into account that the class of fibrations in a cofibrantly generated model category coincides with the class $J$-inj, see Definition 2.1.17 (3) in \cite{Hovey1}, we get that
  $$
  \hbox{$J^+_T$-inj}=\{ f:X\to Y \; \, \hbox{in} \; \, \bcS \mid \forall n>0 \; \, \Ev _n(f)\; \, \hbox{is a fibration in} \; \bcC \} \; ,
  $$
i.e. the class $J^+_T$-inj is the class of positive level fibrations in $\bcS $.

Similarly,
  $$
  \hbox{$I^+_T$-inj}=\{ f:X\to Y \; \, \hbox{in} \; \, \bcS \mid \forall n>0 \; \, \Ev _n(f)\; \, \hbox{is a trivial fibration in}\, \bcC \} \; .
  $$

It follows that
  $$
  \hbox{$I^+_T$-inj} \subset \hbox{$J^+_T$-inj} \; .
  $$
By definition, it means that all morphisms in $I^+_T$-inj have the right lifting property with respect to all morphisms from $J_T^+$. Then it means that
  $$
  J^+_T\subset \hbox{$I^+_T$-cof} \; ,
  $$
as required.

As a result,
  $$
  \hbox{$J^+_T$-cell} \subset W_T^+\cap \hbox{$I_T^+$-cof} \; ,
  $$
and the fourth condition is done.

\bigskip

\noindent {\it Fifth and sixth condition}

\medskip

\noindent The above descriptions of the classes $J_T^+$-inj and $I_T^+$-inj give that
  $$
  \hbox{$J_T^+$-inj} \, \cap \, W_T^+ = \hbox{$I_T^+$-inj } \; .
  $$
This gives the conditions five and six in Theorem 2.1.19 in Hovey's book.

\medskip

Thus, the sets $I_T^+$, $J_T^+$ and $W_T^+$ generate a model structure in $\bcS $, denoted by $\bcM ^+$, such that weak equivalences in it are those morphisms $f:X\to Y$ in which $f_n:X_n\to Y_n$ is a weak equivalence in $\bcC $ for all $n>0$.

\end{pf}

\begin{corollary}
\label{fibcofibplus}
A morphism $f:X\to Y$ in $\bcS $ is a fibration in $\bcM ^+$ if and only if $f_n:X_n\to Y_n$ is a fibration in $\bcC $ for any $n>0$. A morphism $f:X\to Y$ in $\bcS $ is a cofibration in $\bcM ^+$ if and only if $f$ is a cofibration in $\bcM $ and $f_0:X_0\to Y_0$ is an isomorphism. In particular, an object $X$ in $\bcS $ is cofibrant in $\bcM ^+$ if and only if $X$ is cofibrant in $\bcM $ and $X_0=*$.
\end{corollary}

\begin{pf}
The corollary can be proved using the definition of $I_T^+$, $J_T^+$, left lifting and the adjunction between $F_n$ and $\Ev _n$.
\end{pf}

\bigskip

\section{Loop spectra}
\label{loops}

Let $\bcD $ be a simplicial closed symmetric monoidal model category. In particular, for any object $X$ in $\bcD $ the functor $-\wedge X$ has right adjoint functor $\iHom (X,-)$. This is nothing but the function object whose value $\iHom (X,Y)$, for any object $Y$ in $\bcD $, can be viewed as ``functions" from $X$ to $Y$. Certainly, $\iHom (-,-)$ is a bifunctor from $\bcD ^{\op }\times \bcD $ to $\bcD $.

Being a simplicial category, $\bcD $ also has a bifunctor $\Map (-,-)$ from $\bcD ^{\op }\times \bcD $ to the category of simplicial sets $\SSets $ with all nice adjunctions, see \cite{Hovey1} and \cite{GJ}. Since the setting is symmetric and simplicial, we will systematically ignore the difference between the left and right versions of $\iHom $ and $\Map $, see a remark on page 131 in \cite{Hovey1}.

For any simplicial set $U$ we have that its $n$-slice $U_n$ is
canonically isomorphic to the Hom-set $\Hom _{\SSets }(\Delta
^n,U)$. Using the adjunction between $\Map (X,-)$ and $X\wedge -$,
see \cite{GJ}, we obtain that $\Hom _{\SSets }(\Delta ^n,\Map
(X,Y))$ is isomorphic to $\Hom _{\bcD }(X\wedge \Delta ^n,Y)$. Then,
  $$
  \Map (X,Y)_n\simeq \Hom _{\bcD }(X\wedge \Delta ^n,Y)\; .
  $$

Objects $\Map (X,Y)$ come from the simplicial structure of the
category $\bcD $. To provide them with a homotopical meaning we need
to replace $X$ and $Y$ by their cofibrant and fibrant replacements
$QX$ and $RY$ respectively. Then let
  $$
  \map (X,Y) = \Map (QX,RY)\; ,
  $$
so that we obtain yet another bifunctor $\map (-,-)$ from the category $\bcD ^{\op }\times \bcD $ to $\SSets $, see \cite{Hovey2}, Section 2.

Now let $\bcD $ be the category of symmetric spectra $\bcS $. Let $Q$ and $R$ be the cofibrant and fibrant replacement functors with respect to the model structure $\bcM $, and let $Q^+$ and $R^+$ be the cofibrant and fibrant replacement functors with respect to the model structure $\bcM ^+$. Cofibrations do not change when passing to localizations, so that $Q$ remains the same in the localizations of the model structure $\bcM $ by $S$ or $S^+$, and $Q^+$ remains the same in the localizations of the model structure $\bcM ^+$ by $S$ or $S^+$. Respectively, we define two bifunctors
  $$
  \map (X,Y) = \Map (QX,RX)
  $$
and
  $$
  \map ^+(X,Y) = \Map (Q^+X,R^+X)
  $$
from $\bcS ^{\op }\times \bcS $ to $\SSets $.

\bigskip

Next, following \cite{HSS} (and \cite{MMSS}), for any spectrum $X$ in $\bcS $ let
  $$
  \Theta X := \iHom (F_1(T),X)\; ,
  $$
and let
  $$
  \theta :X\lra \Theta X
  $$
be a morphism induced by the morphism $\zeta ^{\uno }_0:F_1(T)\to F_0(\uno )$.

It is useful to interpret the functor $\Theta $ as a loop spectrum. Indeed, if
  $$
  s_- : \bcS \lra \bcS
  $$
is a shift functor
  $$
  s_-=\iHom (F_1(\uno ),-)\; ,
  $$
see Definition 8.9 in \cite{Hovey2}, then $\Theta $ is isomorphic to the composition of $s_-$ and a loop-spectrum functor
  $$
  (-)^T=\iHom (F_0(T),-) : \bcS \lra \bcS \; ,
  $$
loc.cit.

We also have iterations
  $$
  \Theta ^0X=X\; ,
  $$
  $$
  \Theta ^nX := \Theta (\Theta ^{n-1}X)\; ,
  $$
and
  $$
  \theta ^n : X\lra \Theta ^nX\; ,
  $$
being a composition of morphisms $\Theta ^i(\theta ):\Theta ^iX\to \Theta ^{i+1}X$ for all $i=0,\dots ,n-1$.

We can also take the colimit
  $$
  \Theta ^{\infty }X=\colim _n\, \Theta ^nX
  $$
with respect to the morphisms $\Theta ^i(\theta )$, and consider the corresponding morphism
  $$
  \theta ^{\infty } : X\lra \Theta ^{\infty }X\; .
  $$

\bigskip

The meaning of the above constructions comes from topology.
Indeed, let $\bcC $ be the category of pointed simplicial sets $\SSets _*$, $\uno $ be the colon $S^0$, $T$ be the simplicial circle $\Delta [1]/\partial \Delta [1]$, so that $\bcS $ is the category of topological symmetric spectra from \cite{HSS}. For any pointed simplicial set $Y$ let
  $$
  X=F_0(Y)=\Sigma ^{\infty }Y
  $$
be the symmetric $S^1$-suspension spectrum of $Y$. Then
  $$
  \Map (F_1(S^1),X)\simeq \Map (S^1,\Ev _1X)=
  $$
  $$
  =\Map (S^1,S^1\wedge Y)=\Omega \Sigma Y
  $$
-- the simplicial set of loops in the suspension $\Sigma Y$ of the pointed simplicial set $Y$. By adjunction between $F_0$ and $\Ev _0$ we have that
  $$
  Y\simeq \Map (S^0,Y)\simeq \Map (F_0(S^0),F_0(Y))=
  \Map (F_0(S^0),X)\; .
  $$
As the suspension $\Sigma $ is left adjoint to the loop-functor $\Omega $, the identity morphism $\id :\Sigma Y\to \Sigma Y$ gives a morphism $\theta ':Y\to \Omega \Sigma Y$. In view if the above isomorphisms, $\theta '$ is nothing but the morphism
$\Map (\zeta _0^{S^0},X)$, induced by the morphism $\zeta _0^{S^0}:F_1(S^1)\to F_0(S^0)$. In other words, $\theta $ is a ``spectralized" morphism $\theta '$ obtained by replacing Homs by internal Homs in $\bcS $.

Iterating the process we would see that the morphisms $\theta ^n:X\to \Theta ^nX$ come from the morphisms $Y\to \Omega ^n\Sigma ^nY$, and the morphism $\theta ^{\infty }:X\to \Theta ^{\infty }X$ comes from the morphism $Y\to \Omega ^{\infty }\Sigma ^{\infty }Y$ in topology, where the simplicial set $\Omega ^{\infty }\Sigma ^{\infty }Y$ is sometimes denoted by $QY$.

\bigskip

If we will do the same construction $\theta ^{\infty }:X\to \Theta ^{\infty }X$ in the category of non-symmetric spectra over $\SSets _*$, then, as far as we can see, $\Theta ^{\infty }X$ will be an $\Omega $-spectrum, and $\Theta ^{\infty }(-)$ will be a fibrant replacement functor for non-symmetric spectra over $\SSets _*$, see \cite{BF}. However, in symmetric spectra over $\SSets _*$, $\Theta ^{\infty }X$ need not be an $\Omega $-spectrum, and $\theta ^{\infty }$ need not be a stable equivalence, \cite{HSS}.

\bigskip

Now we come back to the category $\bcS $ of abstract symmetric $T$-spectra over $\bcC $.

\begin{proposition}
\label{lemma5}
Let $X$ be an $S^+$-local object in $\bcS $ with respect to the positive projective model structure $\bcM ^+$. Then:

\begin{itemize}

\item[]{(i)}
$\Theta X$ is an $S$-local object with respect to the projective model structure $\bcM $, and

\item[]{(ii)}
the morphism $\theta :X\to \Theta X$ is a weak equivalence in the model structure $\bcM ^+$.
\end{itemize}

\end{proposition}

\begin{pf}
First of all we need to show that $\Theta X$ is fibrant in $\bcM $. Let $f:A\to B$ be a trivial cofibration in $\bcM $, and consider the following commutative square
  $$
  \diagram
  A \ar[dd]_-{f} \ar[rr]^-{} & & \Theta X \ar[dd]^-{} \\ \\
  B \ar@{.>}[uurr]^-{\exists h}\ar[rr]^-{} & & \ast
  \enddiagram
  $$
We need to find a morphism $h:B\to \Theta X$ completing the diagram to be commutative. By adjunction between $-\wedge F_1(T)$ and $\iHom (F_1(T),-)$ the lifting $h$ exists if and only if there exists a lifting $h'$ making the diagram
  $$
  \diagram
  A\wedge F_1(T) \ar[dd]_-{f\wedge F_1(T)} \ar[rr]^-{} & & X \ar[dd]^-{} \\ \\
  B\wedge F_1(T) \ar@{.>}[uurr]^-{\exists h'}\ar[rr]^-{} & & \ast
  \enddiagram
  $$
commutative. The object $F_1(T)$ is cofibrant in $\bcM $ because $T$ is cofibrant and the functor $F_1$ is left Quillen with respect to the model structure $\bcM $. Then $f\wedge F_1(T)$ is a trivial cofibration in $\bcM $. The specificity of the spectrum $F_1(T)$ yields that $(A\wedge F_1(T))_0=\ast $ and $(B\wedge F_1(T))_0=\ast $, so that $(f\wedge F_1(T))_0$ is an isomorphism. Therefore, $f\wedge F_1(T)$ is a trivial cofibration not only in $\bcM $ but also in $\bcM ^+$. Since $X$ is fibrant in $\bcM ^+$, because it is $S^+$-local with respect to $\bcM ^+$ by assumption, the required $h'$ exists.

Thus, $\Theta X$ is fibrant in $\bcM $, and we can start to prove the first part of the proposition. In order to show that $\Theta X$ is $S$-local, with respect to $\bcM $, we need to show that for any cofibrant object in $U$ in $\bcM $, and for any non-negative integer $n$ the morphism
  $$
  (\zeta _n^U)^* : \map (F_{n+1}(T\wedge U),\Theta X)\lra
  \map (F_n(U),\Theta X)
  $$
is a weak equivalence of simplicial sets. As $U$ is cofibrant, the spectrum $F_n(U)$ is cofibrant in $\bcM $ too. The spectrum $\Theta X$ is fibrant in $\bcM $. Therefore, the simplicial set
  $$
  \map (F_n(U),\Theta X)
  $$
is weak equivalent to the simplicial set
  $$
  \Map (F_n(U),\Theta X)\simeq
  \Hom _{\bcS }(F_n(U)\wedge \Delta ^{\bullet },\Theta X)
  $$
By the adjunction between $-\wedge F_1(T)$ and $\iHom (F_1(T),-)$ we get an isomorphism
  $$
  \Hom _{\bcS }(F_n(U)\wedge \Delta ^{\bullet },\Theta X)\simeq
  \Hom _{\bcS }(F_n(U)\wedge F_1(T)\wedge
  \Delta ^{\bullet },X)\; .
  $$
But
  $$
  \Hom _{\bcS }(F_n(U)\wedge F_1(T)\wedge
  \Delta ^{\bullet },X)\simeq \Map (F_n(U)\wedge F_1(T),X)\; .
  $$
Besides, the spectrum $F_n(U)\wedge F_1(T)$ is cofibrant in $\bcM ^+$, and $X$ is fibrant in $\bcM ^+$ by assumption, so that
  $$
  \Map (F_n(U)\wedge F_1(T),X)\sim
  \map ^+(F_n(U)\wedge F_1(T),X)\; ,
  $$
where $\sim $ stays for weak equivalences of simplicial sets. As a result, we obtain that
  $$
  \map (F_n(U),\Theta X)\sim
  \map ^+(F_n(U)\wedge F_1(T),X)\; .
  $$

Similarly, we get an isomorphism
  $$
  \map (F_{n+1}(U\wedge T),\Theta X)\sim \map ^+(F_n(U\wedge T)\wedge F_1(T),X)\; .
  $$

Consider now a commutative square
  $$
  \diagram
  \map (F_n(U),\Theta X) \ar[dd]_-{\sim } \ar[rr]^-{(\zeta _n^U)^*} & & \map (F_{n+1}(U\wedge T),\Theta X)
  \ar[dd]^-{\sim } \\ \\
  \map ^+(F_n(U)\wedge F_1(T),X) \ar[rr]^-{(\zeta _{n+1}^{U\wedge T})^*} & & \map ^+(F_{n+1}(U\wedge T)\wedge F_1(T),X)
  \enddiagram
  $$
As $X$ is $S^+$-local with respect to $\bcM ^+$ by assumption, and the morphism $\zeta _{n+1}^{U\wedge T}$ is in $S^+$, the bottom horizontal morphism is the above diagram is a weak equivalence of simplicial sets. And, as we have seen just now, the vertical morphisms are weak equivalences of simplicial sets. Then the top horizontal morphism $(\zeta _n^U)^*$ is a weak equivalence of simplicial sets as well, and $(i)$ is done.

To prove $(ii)$ all we need to show is that the morphism $\theta _n:X_n\to (\Theta X)_n$ is a weak equivalence in $\bcC $ for all $n\geq 1$. Since $X$ is fibrant in $\bcM ^+$ by assumption, $X_n$ is fibrant in $\bcC $, provided $n\geq 1$. The object $(\Theta X)_n$ is fibrant in $\bcC $ because $\Theta X$ is fibrant in $\bcM $ by $(i)$. Therefore, it is enough to show that for any cofibrant object $B$ in $\bcC $ the corresponding morphism
  $$
  (\theta _n)_* : \map (B,X_n)\lra \map (B,(\Theta X)_n)
  $$
is a weak equivalence of simplicial sets.

Recall that $\bcC $ is simplicial. As $B$ is cofibrant and $\Ev _n(X)$ is fibrant in $\bcC $, we have that
  $$
  \map (B,\Ev _n(X))\sim \Map (B,\Ev _n(X))=\Hom _{\bcC }(B\wedge \Delta ^{\bullet },\Ev _n(X))\; .
  $$
By the adjunction between $F_n$ and $\Ev _n$ we have that
  $$
  \Hom _{\bcC }(B\wedge \Delta ^{\bullet },\Ev _n(X))\simeq
  \Hom _{\bcS }(F_n(B\wedge \Delta ^{\bullet }),X)
  $$
But
  $$
  F_n(B\wedge \Delta ^{\bullet })\simeq F_n(B)\wedge \Delta ^{\bullet }
  $$
by the definition of the action of simplicial sets on spectra. Therefore, we obtain
  $$
  \map (B,\Ev _n(X))\sim \Hom _{\bcS }(F_n(B\wedge \Delta ^{\bullet }),X)\simeq \Hom _{\bcS }(F_n(B)\wedge \Delta ^{\bullet },X)=
  $$
  $$
  =\Map (F_n(B),X)\; .
  $$

Similarly, we get a weak equivalence
  $$
  \map (B,\Ev _n(\Theta X))\sim \Map (F_n(B),\Theta X)\; .
  $$

Now we look at the commutative square
  $$
  \diagram
  \map (B,\Ev _n(X)) \ar[dd]_-{\sim } \ar[rr]^-{(\theta _n)_*} & & \map (B,\Ev _n(\Theta X)) \ar[dd]^-{\sim } \\ \\
  \Map (F_n(B),X) \ar[rr]^-{\theta _*} & &
  \Map (F_n(B),\Theta X)
  \enddiagram
  $$
Changing the bottom horizontal row by means of the adjunction
  $$
  \Map (F_n(B),\Theta X)\simeq \Map (F_n(B)\wedge F_1(T),X)\simeq \Map (F_{n+1}(B\wedge T),X)
  $$
we get a new commutative square
  $$
  \diagram
  \map (B,X_n) \ar[dd]_-{\sim } \ar[rr]^-{(\theta _n)_*} & & \map (B, \Theta X_n) \ar[dd]^-{\sim } \\ \\
  \Map (F_n(B),X) \ar[rr]^-{(\zeta _n^B)^*} & & \Map (F_{n+1}(B\wedge T),X)
  \enddiagram
  $$
As $n\geq 1$, the morphism $\zeta _n^B$ is in $S^+$, the objects $F_n(B)$ and $F_{n+1}(B\wedge T)$ are cofibrant in $\bcM ^+$. Besides, $X$ is $S^+$-local. Therefore, the bottom horizontal morphism in the last commutative square is a weak equivalence of simplicial sets. Since the vertical morphisms are weak equivalences, we obtain that the top horizontal $(\theta _n)_*$ is a weak equivalence of simplicial sets.
\end{pf}

If $\bcA $ and $\bcA '$ are two model structures in the same category $\bcB $ then we will use the symbols $Ho(\bcA )$ and $Ho(\bcA ')$ for the homotopy categories of the category $\bcB $ with respect to the model structures $\bcA $ and $\bcA '$ respectively. We also will be using the following lemma.

\begin{lemma}
\label{barsuk}
The functor $\Theta :\bcS \to \bcS $ carries weak equivalences between fibrant objects in $\bcM ^+$ in to weak equivalences in $\bcM $, so that the pair of functors
  $$
  (-\wedge F_1(T),\Theta )
  $$
is a Quillen adjunction between $\bcM $ and $\bcM ^+$. In particular, there exists right derived functor $R\Theta :Ho (\bcM ^+)\to Ho(\bcM )$.
\end{lemma}

\begin{pf}
Let $f$ be a (trivial) cofibration in the model structure $\bcM $.
As the model structure $\bcM $ is compatible with the monoidal
structure in $\bcS $, the morphism $f\wedge F_1(T)$ is also a
(trivial) cofibration in $\bcM $. Since $(F_1(T))_0=*$ the morphism
$(f\wedge F_1(T))_0$ is an isomorphism. Since $-\wedge F_1(T)$ has
right adjoint $\Theta $, we are done.
\end{pf}

\section{Positive weak equivalences are stable}
\label{stablepositive}

In this section we will show that any weak equivalence in the positive model structure is a weak equivalence in the stable model structure. This is a consequence of the previous results and the following general effect.

\begin{lemma}
\label{enot}
Let $\bcD $ be a closed symmetric monoidal model category with a product $\wedge $ and unit $\uno $. Suppose $\bcD $ is cofibrantly generated and that the domains of the generating cofibrations are cofibrant. Let $U$ be cofibrant, $X$, $Y$ fibrant objects, and let
  $$
  u : U\to \uno\; ,\quad f:X\to Y
  $$
be two morphisms, all in $\bcD $. Denote by $(u)$ the set of morphisms
  $$
  V\wedge u : V\wedge U\to V\; ,
  $$
where $V$ runs through domains and codomains of generating cofibrations in $\bcD $. Suppose furthermore that the morphism
  $$
  f_* : \iHom (U,X)\lra \iHom (U,Y)
  $$
is a weak equivalence in $\bcD $. Then, $f$ is a weak equivalence in the Bousfield localized category $\bcD _{(u)}$.
\end{lemma}

\begin{pf}
The category $\bcD _{(u)}$ is closed monoidal model by Proposition 33 in \cite{GG2}, and for any cofibrant object $V$ in $\bcD $ the morphism $V\wedge u$ is a weak equivalence in $\bcD _{(u)}$ by Lemma 32 in loc.cit.

Let $q:Q\uno \to \uno $ be a cofibrant replacement of the unit in
$\bcD $. In the commutative diagram
  $$
  \diagram
  Q\uno \wedge U\ar[dd]_-{q\wedge U}
  \ar[rr]^-{Q\uno \wedge u} & &
  Q\uno \ar[dd]^-{q} \\ \\
  U \ar[rr]^-{u} & & \uno
  \enddiagram
  $$
the morphism $q$ is a weak equivalence in $\bcD $ by definition, the morphism $q\wedge U$ is a weak equivalence in $\bcD $ by one of the axioms of the monoidal model structure, and $Q\uno \wedge u$ is a weak equivalence in $\bcD_{(u)}$. Therefore, $u$ is a weak equivalence in $\bcD_{(u)}$.

The morphism $u$ defines a morphism
  $$
  u^* : X\simeq \iHom(\uno ,X)\to \iHom(U,X)\; .
  $$
Let
  $$
  r:X\lra R_{(u)}X
  $$
be the fibrant replacement in $\bcD_{(u)}$. As $R_{(u)}X$ is
fibrant, $U$ is cofibrant and $u$ is a weak equivalence in the
closed monoidal model category $\bcD _{(u)}$, the morphism
  $$
  u^* : R_{(u)}X\simeq \iHom (\uno ,R_{(u)}X))\lra \iHom (U,R_{(u)}X)
  $$
is a weak equivalence in $\bcD _{(u)}$ by \cite[Lemma 4.2.7]{Hovey1}. The morphism $r$ is a weak equivalence in $\bcD _{(u)}$ by definition. As the square
  $$
  \diagram
  X \ar[dd]_-{r} \ar[rr]^-{u^*} & &
  \iHom (U,X) \ar[dd]^-{r_*} \\ \\
  R_{(u)}X \ar[rr]^-{u^*} & &
  \iHom (U,R_{(u)}X)
  \enddiagram
  $$
is commutative, the composition
  $$
  X\stackrel{u^*}\lra \iHom (U,X)
  \stackrel{r_*}\lra \iHom (U,R_{(u)}X)
  $$
is an isomorphism in the homotopy category $Ho(\bcD _{(u)})$, which means that $X$ is functorially a retract of
$\iHom (U,X)$ in $Ho(\bcD _{(u)})$.

In particular, $f$ is a retract of an isomorphism $f_*:\iHom (U,X)\to \iHom (U,Y)$ in $Ho(\bcD_{(u)})$. As a retract of an isomorphism is an isomorphism, $f$ is an isomorphism in $Ho(\bcD _{(u)})$, and so it is a weak equivalence in $\bcD _{(u)}$.
\end{pf}

\bigskip

\begin{proposition}
\label{important}
Any positive weak equivalence is a stable weak equivalence.
\end{proposition}

\begin{pf}
Let $f:X\to Y$ be a positive weak equivalence, i.e. $f_n:X_n\to Y_n$ is a weak equivalence in $\bcC $ for any $n>0$. We are going to apply Lemma \ref{enot} when
  $$
  U=F_1(T)
  $$
and
  $$
  u=\zeta _0^{\uno }:F_1(T)\lra F_0(\uno )=\uno \; .
  $$
Notice that $U$ is cofibrant in $\bcM $ and, without loss of
generality, we may assume that $X$ and $Y$ are fibrant objects in $\bcM $, because fibrant replacements in $\bcM $ are level
equivalences and do not change neither the condition of the proposition, nor its conclusion. Then $X$ and $Y$ are fibrant in $\bcM ^+$, too. As $f$ is a weak equivalence in $\bcM ^+$, by Lemma \ref{barsuk}, the morphism $\Theta f=\iHom (F_1(T),f)$ is a weak equivalence in $\bcM $. Then $f$ is a weak equivalence in the model structure $\bcM _{(\zeta _0^{\uno })}$ by Lemma \ref{enot}. To complete the proof we need only to observe that, for any cofibrant object $V$ in $\bcM $, the
morphism $V\wedge \zeta _0^{\uno }$ is a stable weak equivalence, so that $(\zeta _0^{\uno })$ consists of weak equivalences in $\bcM _S$. Actually, $\bcM _{(\zeta ^{\uno }_0)}=\bcM _S$, because  $\zeta ^X_n=F_n(X)\wedge \zeta ^{\uno }_0$.
\end{pf}

Recall that $Q$ is the cofibrant replacement functor with respect to the model structure $\bcM $, and $Q^+$ is the cofibrant replacement functor with respect to the model structure $\bcM ^+$. Replacing $Q^+$ by $Q^+Q$, we obtain a natural transformation
  $$
  Q^+\lra Q\; .
  $$

\begin{corollary}
\label{kozerog}
Let $X$ and $Z$ be two objects in $\bcS $, such that $Z$ is $S$-local with respect to the projective model structure $\bcM $ in $\bcS $. Then the morphism
  $$
  \map (X,Z)\lra \map ^+(X,Z)\; ,
  $$
induced by the above natural morphism $Q^+X\to QX$, is a weak equivalence of simplicial sets.
\end{corollary}

\begin{pf}
As $Z$ is $S$-local with respect to $\bcM $, it is fibrant in $\bcM $, and so in $\bcM ^+$. Let
  $$
  q:QX\lra X
  $$
be the cofibrant replacement in $\bcM $. The morphisms
  $$
  q^* : \map (X,Z)\lra \map (QX,Z)
  $$
and
  $$
  q^* : \map ^+(X,Z)\lra \map ^+(QX,Z)
  $$
are both weak equivalences of simplicial sets. Therefore, without loss of generality, one can assume that $X$ is cofibrant in $\bcM $.

Let now
  $$
  q^+ : Q^+X\lra X
  $$
be the cofibrant replacement of $X$ in $\bcM ^+$. Then $q^+$ is a positive weak equivalence, hence a stable weak equivalence in $\bcM _S$, by Proposition \ref{important}. The objects $X$ and $Q^+X$ are cofibrant in $\bcM $, so in $\bcM _S$, and $Z$ is fibrant in $\bcM _S$. Then the morphism
  $$
  \map (X,Z)\sim \Map (X,Z)\stackrel{(q^+)^*}{\lra }\Map (Q^+X,Z)\sim \map ^+(X,Z)
  $$
is a weak equivalence of simplicial sets because $\bcS $ is a simplicial model category with respect to the model structure $\bcM _S$.
\end{pf}

\begin{remark}
{\rm For a natural $n$ call an $n$-level weak equivalence (fibration) a morphism in $\bcS $ which is a level weak equivalence (fibration) for $i$-slices with $i\geq n$. These two classes of morphisms define a model structure $\bcM ^{\geq n}$ on $\bcS $. Cofibrations in $\bcM ^{\geq n}$ are cofibrations in $\bcM $ which are isomorphisms on $i$-slices with $i<n$ and $n$-level weak equivalences. By methods similar to those used above one show that any $n$-level weak equivalence is a stable weak equivalence. }
\end{remark}

\section{Main theorem}

Recall that $W_{T,S}$ is the set of weak equivalences in $\bcM _S$, and $W^+_{T,S^+}$ is the set of weak equivalences in $\bcM ^+_{S^+}$. Let also $W^+_{T,S}$ be the set of weak equivalences in $\bcM ^+_S$.

\begin{theorem}
\label{main}
In the notation above,
  $$
  W_{T,S}=W^+_{T,S^+}=W^+_{T,S}\; .
  $$
\end{theorem}

\begin{pf}
Let $f:X\to Y$ be a weak equivalence in $\bcM _S$. In order to prove that $f$ is a weak equivalence in $\bcM ^+_{S^+}$ we need to show that for any $S^+$-local object $Z$ in $\bcM ^+$ the morphism
  $$
  \map ^+(Y,Z)\lra \map ^+(X,Z)
  $$
is a weak equivalence of simplicial sets. The morphism
  $$
  \theta : Z\lra \Theta Z\; ,
  $$
together with the morphism $f$, give rise to the commutative square
  $$
  \diagram
  \map ^+(Y,Z) \ar[dd]_-{\theta _*} \ar[rr]^-{f^*} & &
  \map ^+(X,Z) \ar[dd]^-{\theta _*} \\ \\
  \map ^+(Y,\Theta Z) \ar[rr]^-{f^*} & & \map ^+(X,\Theta Z)
  \enddiagram
  $$
As $Z$ is $S^+$-local in $\bcM ^+$, Proposition \ref{lemma5} (i) gives that $\Theta Z$ is $S$-local in $\bcM $. Since $f$ is a weak equivalence in $\bcM _S$, the morphism
  $$
  f^* : \map (Y,\Theta Z)\lra \map (X,\Theta Z)
  $$
is a weak equivalence of simplicial sets. Applying Corollary \ref{kozerog} we obtain that the lower $f^*$ in the above commutative square is also a weak equivalence of simplicial sets. Proposition \ref{lemma5} (ii) gives that the morphism $\theta $ is a weak equivalence in $\bcM ^+$. It follows that the vertical morphisms in the above commutative square are weak equivalences of simplicial sets. Then the top horizontal morphism is a weak equivalence of simplicial sets, as required. Thus, $W_{T,S}\subset W^+_{T,S^+}$.

Let $f:X\to Y$ be a weak equivalence in $\bcM ^+_{S^+}$. We want to show that $f$ is a weak equivalence in $\bcM _S$. Take any $S$-local object $Z$ in $\bcM $ and look at the commutative diagram
  $$
  \diagram
  \map (Y,Z) \ar[dd]_-{} \ar[rr]^-{f^*} & &
  \map (X,Z) \ar[dd]^-{} \\ \\
  \map ^+(Y,Z) \ar[rr]^-{f^*} & & \map ^+(X,Z)
  \enddiagram
  $$
As $Z$ is $S$-local in $\bcM $, it is $S^+$-local in $\bcM ^+$. Since $f$ is a weak equivalence in $\bcM ^+_{S^+}$, the lower horizontal morphism is a weak equivalence of simplicial sets. The vertical arrows in the diagrams are isomorphisms from Corollary \ref{kozerog}. Then the top horizontal arrow is a weak equivalence of simplicial sets, for any $S$-local object $Z$ in $\bcM $. It means that $f$ is a weak equivalence in $\bcM _S$.

Thus, $W_{T,S}=W^+_{T,S^+}$. In particular, all morphisms in $S$ are weak equivalences in $\bcM ^+_{S^+}$. This implies that $(\bcM ^+_{S^+})_S=\bcM ^+_{S^+}$. On the other hand, $(\bcM ^+_{S^+})_S=\bcM ^+_S$, because $S^+\subset S$.
\end{pf}

\bigskip

\bigskip

\begin{small}

\end{small}

\vspace{4mm}

\begin{small}

{\sc Steklov Mathematical Institute, 8 Gubkina Street,
Moscow, 119991 Russia}

\end{small}

\begin{small}
{\sc AG Laboratory, HSE, 7 Vavilova Street, Moscow, 117312 Russia}

\end{small}

\begin{footnotesize}

{\it E-mail address}: {\tt gorchins@mi.ras.ru}

\end{footnotesize}

\bigskip

\begin{small}

{\sc Department of Mathematical Sciences, University of Liverpool,
Peach Street, Liverpool L69 7ZL, England, UK}

\end{small}

\begin{footnotesize}

{\it E-mail address}: {\tt vladimir.guletskii@liverpool.ac.uk}

\end{footnotesize}

\end{document}